\documentclass{ifacconf}
\usepackage{amssymb,amsmath,amsfonts}
\usepackage{mathtools}
\usepackage{graphicx}
\graphicspath{{Images/}}
\usepackage{tikz}
\usepackage{pgfplots}
\pgfplotsset{compat=1.14}
\usepackage[utf8]{inputenc}
\usepackage{microtype}
\usepackage{nicefrac}
\usepackage{enumitem}
\usepackage{xspace}
\usepackage{url}
\usepackage{natbib}
\setcitestyle{round}
\usepackage[absolute]{textpos}

\newcommand{\eg}{\emph{e.g.,}\xspace}
\newcommand{\ie}{\emph{i.e.,}\xspace}

\newcommand{\tradeoff}{communication-computation trade-off\xspace}
\newcommand{\processing}{processing\xspace}
\newcommand{\preprocessing}{preprocessing\xspace}

\newcommand{\Real}[1]{ { {\mathbb R}^{#1} } }
\newcommand{\Realp}[1]{ { {\mathbb R}^{#1}_+ } }

\DeclareMathOperator*{\argmin}{arg\,min}
\newcommand{\psteady}{p_{\infty|\infty-\tau_{\textit{tot}}}}
\newcommand{\psteadytau}{p_{\infty|\infty-\tau}}
\newcommand{\psteadytautauc}{p_{\infty|\infty-{\tau_{s}}}}

\newtheorem{ass}{Assumption}

\newcommand{\blue}[1]{{\color{blue}#1}}

\newcommand{\linkToPdf}[1]{\href{#1}{\blue{(pdf)}}}
\newcommand{\linkToPpt}[1]{\href{#1}{\blue{(ppt)}}}
\newcommand{\linkToCode}[1]{\href{#1}{\blue{(code)}}}
\newcommand{\linkToWeb}[1]{\href{#1}{\blue{(web)}}}
\newcommand{\linkToVideo}[1]{\href{#1}{\blue{(video)}}}
\newcommand{\linkToMedia}[1]{\href{#1}{\blue{(media)}}}
\newcommand{\award}[1]{\xspace} 

\begin{document}
	\begin{frontmatter}
		
		\begin{textblock}{20}(-2,1.2)
			\centering
			This paper has been accepted for the 21st IFAC World Congress in Berlin, July 12-17, 2020.\\
			Please cite the paper as: L. Ballotta, L. Schenato, and L. Carlone,\\
			“From Sensor to Processing Networks: Optimal Estimation with Computation and Communication Latency”,\\
			IFAC 2020 World Congress, 2020
		\end{textblock}
		
		\title{From Sensor to Processing Networks: Optimal Estimation with Computation and Communication Latency\thanksref{ack}}
		
		\thanks[ack]{This work was partially funded by the ONR RAIDER program (N00014-18-1-2828), the CARIPARO Foundation Visiting Programme ``HiPeR" and the Italian Ministry of Education PRIN n. 2017NS9FEY.}
		
		\author[unipd]{Luca Ballotta}
		\author[unipd]{Luca Schenato} 
		\author[MIT]{Luca Carlone} 
		
		\address[unipd]{Department of Information Engineering, University of Padova, Padova, 35131, Italy (e-mail: \{ballotta, 	schenato\}@dei.unipd.it)}
		\address[MIT]{Laboratory for Information \& Decision Systems, Massachusetts Institute of Technology, Boston, 02139, USA (e-mail: lcarlone@mit.edu)}


\begin{abstract}
This paper investigates the use of a networked system
(\eg swarm of robots, smart grid, sensor network) to monitor a time-varying phenomenon of interest in the 
presence of communication and computation latency.  
Recent advances in edge computing have enabled processing to be spread across the network, 
hence we investigate the fundamental \tradeoff,
arising when a sensor has to decide whether to transmit raw data 
(incurring communication delay) or preprocess them (incurring computational delay) in order to compute an accurate 
estimate of the state of the phenomenon of interest.
We propose two key contributions. First, we formalize the notion of \emph{processing network}.
Contrarily to \emph{sensor and communication networks}, where the designer is concerned with 
the design of a suitable communication policy, in a \emph{processing network} 
one can also control when and where the computation occurs in the network. 
The second contribution is to provide analytical results on the optimal preprocessing delay (\ie the optimal time spent on computations at each sensor) for the case with a single sensor and  multiple homogeneous sensors.
Numerical results substantiate our claims that accounting for computation latencies (both at sensor and estimator side) and communication delays can largely impact the estimation accuracy.

\end{abstract}

		
	\end{frontmatter}
	

\section{Introduction}

Networked systems are becoming an ubiquitous technology across many application domains, 
including city-wide air-pollution monitoring~\citep{8405565}, smart power grids~\citep{Pasqualetti11cdc-safeCPS}, swarms of mobile robots for target tracking~\citep{4937860}, interconnected autonomous vehicles and self-driving cars~\citep{Shalev-Shwartz17arxiv-safeDriving}.
 Progress on communication systems, such as the development of 5G, carries the promise of further expanding the reach of these systems by enabling more effective and larger-scale deployments.
At the same time, recent advances on embedded computing, 
from embedded GPU-CPU systems to specialized hardware~\citep{Suleiman18jssc-navion}, are now providing unprecedented opportunities for \textit{edge-computing}, where sensor data are processed locally at the sensor to minimize the communication burden.

\begin{figure}
	\centering
	\includegraphics[width=0.9\linewidth]{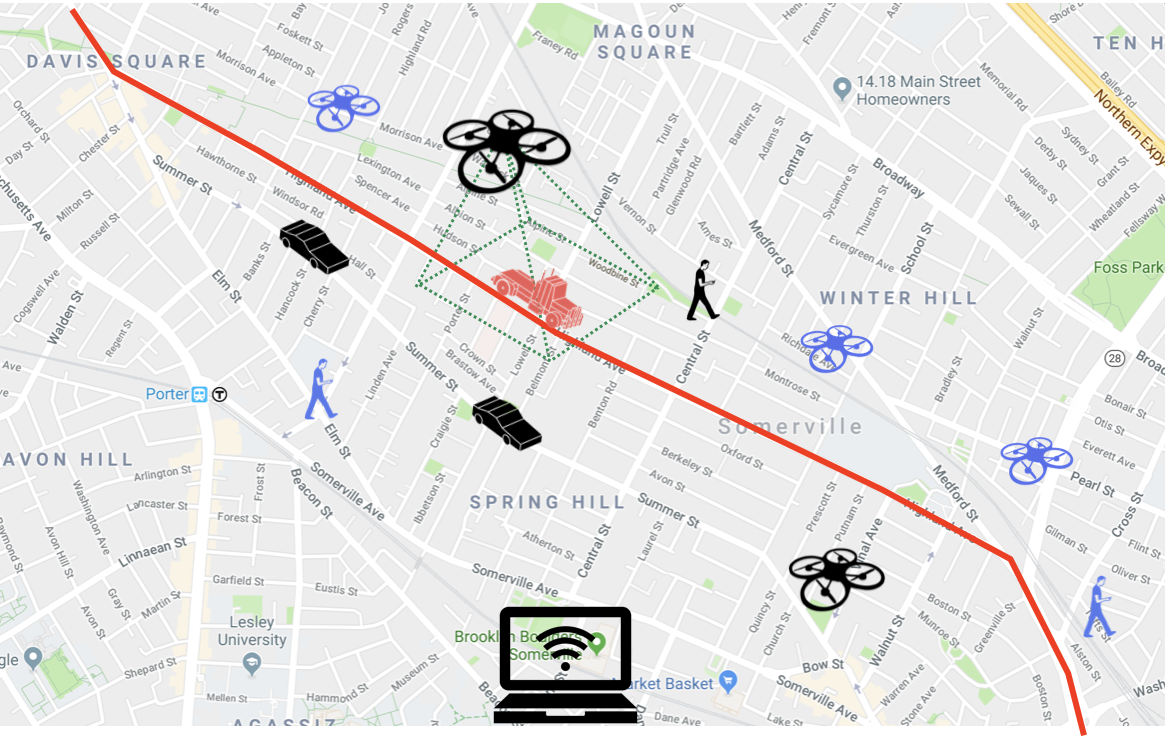}
	\caption{Example of \emph{processing network}: 
	smart sensors (in blue and black) collect, process, and communicate data to
	track the state of a vehicle (in red) in the presence of communication and computation latency.}
	\label{fig:vehicle-tracking}
\end{figure}

The availability of powerful embedded computers creates a nontrivial \textit{\tradeoff}: 
is it best to transmit raw sensor data and incur larger communication and data fusion delays at a central station,
or to perform more preprocessing at the sensors and transmit more refined (less noisy and more compressed) estimates? 
Fig.~\ref{fig:vehicle-tracking} provides an example of this trade-off: 
the figure depicts a network of \emph{smart sensors} (in black and blue) observing and tracking the state of a
moving vehicle (in red) and transmitting data to a central fusion station (the computer at the bottom of the figure), which is in charge of monitoring the state of the red truck.
The smart sensors may have heterogeneous computational resources: for instance, the large drone (in black) might have a powerful onboard GPU-CPU system, while other smart sensors (in blue, e.g., smaller drones, mobile phones) might have limited computation.
Therefore, some sensors might prefer sending raw data and incur larger delays for transmission, while some other sensors might prefer preprocessing the data at the edge. These choices will impact the quality of the red vehicle estimate: larger computation and communication delays will lead to less accurate estimates, hindering the tracking task.

In this paper we investigate the \tradeoff that arises in a networked system 
responsible for estimating the state of a time-varying phenomenon of interest in the presence of computational and communication delays. 
Related work in the IoT community focuses on optimizing data transmission by means of smart communication policies, with respect to estimation performance~\citep{2018arXiv180405618W} or the so-called \emph{Age of Information} (AoI) due to delays and unreliability~\citep{8778671,8469047}.
~\cite{8006543} and~\cite{Bisdikian:2013:QVI:2489253.2489265} introduce the \emph{Value of Information of Update} (VoIU), which addresses the impact new samples have on the current state estimate. Also, the former formalize the \emph{Cost of Update Delays} (CoUD), a non-linear function of AoI expanding such concept, which concurs in the VoIU of samples.
Contrary to this line of work, we focus on monitoring a dynamical system and advocate a unified 
task-driven framework, where computation and communication are jointly modeled in an optimal estimation framework.
In hindsight, we propose a paradigm shift from \emph{sensor and communication networks}, in which one has to decide the best communication policy, to \emph{processing networks}, where one also controls when and where the computation occurs. 
Moreover, we analyze the relation between computation/communication delays and system dynamics, while 
previous work mostly focuses on the channel properties.

Related work in {control, cyber-physical systems, and robotics} focuses on either the co-design of estimation and control in the presence of communication constraints~\citep{Borkar97cccsp-limitedCommControl,Shafieepoorfard13cdc-attentionLQG},
or on the design of the system's sensing and actuation~\citep{Carlone18tro-attentionVIN,Summers16tcns-sensorScheduling}. 
\cite{Tzoumas18acc-sLQG} establish a more direct connection between sensing and estimation performance, by proposing co-design approaches for sensing, estimation, and control.
While these works focus on communication constraints, we attempt to explicitly model 
\emph{computation delays} and understand how they impact the performance of the estimation task.
In robotics, \cite{2019arXiv190205703C} adopt a learning approach for
computational-offloading in cloud-robotics applications.
\cite{Tsiatsis2005} seek a policy to tackle edge-computing delays within a static framework.
\cite{8757960} characterize the performance of resource-constrained devices with cloud fog offloading (with case study on Fast Fourier Transform computation), while~\cite{7800393} investigate multimedia data processing within pipeline and parallel architectures.
Contrary to these works, we consider the system dynamics, we explicitly model communication and computational delays, and we are concerned with the analytical derivation of optimal computation policies for estimation.

We propose the following contributions. First, we formalize the notion of processing network and provide a model which is amenable for analysis (Section~\ref{sec:set-up}). The networked system is modeled as a set of smart sensors in charge of estimating the state of a dynamical system in the presence of communication and computation latency. 
We assume that edge devices run so-called \emph{anytime algorithms}, \ie the quality of their estimates improves with longer runtime. The key idea is to capture the impact of the \preprocessing at each sensor using a processing-dependent measurement noise, such that more \processing leads to more refined measurements. 
Second, we derive fundamental limits for such model: we prove that in two instantiations of the model there is an optimal choice for the amount of preprocessing done at each sensor which can be computed analytically. 
In particular, Sections~\ref{sec:preprocessing-delay}--\ref{sec:preprocessing-comm-delay} consider the continuous-time case with a single sensor and provide closed-form expressions for the optimal computational delay, while Section~\ref{sec:multiple-sensors} generalizes the setup to multiple homogeneous sensors. 
A discussion of potential extensions to heterogeneous sensors and discrete-time systems is briefly presented in Section~\ref{sec:discussion}, while we refer the interested reader to the 
preprint~\cite{2019arXiv191105859B} for a more comprehensive discussion.
Conclusions are drawn in Section~\ref{sec:conclusions}.

\begin{figure*}
	\centering
	\includegraphics[width=0.7\linewidth]{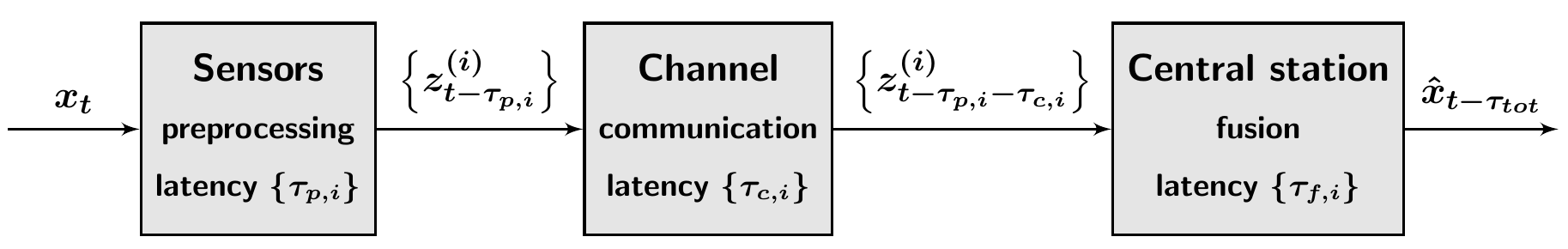}
	\caption{Block diagram of the processing network with latency contributions by preprocessing, communication, and fusion.}
	\label{fig:blocks-delays}
\end{figure*}

\section{Estimation in processing networks: problem formulation}
\label{sec:set-up}

A \emph{processing network} is a set of interconnected \emph{smart sensors} that collect
sensor data and leverage onboard computation to locally preprocess the data before communicating it to a central fusion center.
The goal of the network is to obtain an accurate estimate of the state of a time-varying phenomenon observed by the sensors, in the face of communication and computation latencies.

\subsection{Anatomy of a Processing Network}

\begin{description}[leftmargin=0cm]
	\item[Dynamical system:]
	We consider a processing network monitoring a time-varying phenomenon 
	described by the following linear time-invariant (LTI) stochastic model:
	\begin{equation}
	dx_t = a \, x_tdt + dw_t \label{eq:processModel}
	\end{equation} 
	where $x_t \in \Real{}$ is the to-be-estimated state of the system at time $t$, 
	$a \in \Real{}$ is a constant describing the system dynamics, and $w_t \in \Real{}$ represents process noise. 
	We focus on the scalar system~\eqref{eq:processModel} which can be analyzed analytically and postpone the discussion 
	on the multi-variate case to Section~\ref{sec:multidim}.

	\item[Smart sensors:] 
	The processing network includes $N$ smart sensors,  $\mathcal{N}=\{1,...,N\}$. 
	After acquiring data, each sensor may refine raw measurements via some local preprocessing. 
	For instance, in the robotics application of Fig.~\ref{fig:vehicle-tracking}, each robot is a smart sensor that may process raw data (e.g., images) to obtain 
	local measurements of the state (e.g., the tracked vehicle location in Fig.~\ref{fig:vehicle-tracking}).
	Depending on the time and computational resources, the robot may use more sophisticated algorithms (or a larger number of visual features~\citep{Hartley04}) to obtain more accurate measurements. 
	More generally, the use of \emph{anytime algorithms}~\citep{Zilberstein96ai-anytimeAlgorithms} at each sensor entails a  trade-off, where 
	the more time is spent on preprocessing, the more accurate is the measurement by the sensor. 
	We capture the dependence of the preprocessing time on the refined measurements through the following model:
	\begin{equation}
	\hspace{-1mm}z_t(\tau_p) \!=\! Cx_t\!+\!v_t(\tau_p), \ z_t(\tau_p) \!=\!
	\begin{bmatrix}
	z_t^{(1)}(\tau_{p,1}) & \!... & \!z_t^{(N)}(\tau_{p,N})
	\end{bmatrix}^T \!\label{eq:measurementModel}
	\end{equation}
	where $z_t^{(i)}$ is the measurement collected at time $t$ by the $i$-th sensor,
	$\tau_{p,i}$ is the \textit{preprocessing delay} associated with the $i$-th sensor,  
	and $v_{t}$ is white noise; $\tau_p \doteq \{\tau_{p,i}\}_{i\in\mathcal{N}}$, and $z_t$ contain the delays and measurements from all sensors. 
	In order to capture the anytime nature of the sensor preprocessing, we model the intensity 
	of the white noise $v_{t}$ as a decreasing function of $\tau_p$, see Sections~\ref{sec:preprocessing-delay} and \ref{sec:other-functions}.  		

	\item[Communication:] The sensors send preprocessed data  to the central station for data fusion. To simplify the mathematical analysis, we assume what follows.
	\begin{ass}[Reliable channel]
		Packet loss and channel erasure probabilities are equal to zero.
		\label{ass:model}
	\end{ass}
	This assumption is quite strong in practice, but it is needed for a tractable analytical approach.
	Future work may include more realistic communication models.\\
	Given limited bandwidth, also data transmission induces a \textit{communication delay} $\tau_{c,i}$ for each sensor $i$. In addition to Assumption~\ref{ass:model}, we assume unconstrained channel capacity, so that all sensors can transmit 
	``in parallel''. 
	We then consider two models for $\tau_{c,i}$ as a function of $\tau_{p,i}$:
	\begin{itemize}
		\item \emph{constant $\tau_{c,i}$:} the transmitted packet length/number is fixed and does not depend on the preprocessing;
		in this case the communication delay is a constant, irrespective of the preprocessing delay $\tau_{p,i}$.
		\item \emph{decreasing $\tau_{c,i}$:} in this case, sensor preprocessing \emph{compresses} the measurements, such that a longer preprocessing leads to less packets to transmit.
	\end{itemize}
 	These models are used in Sections~\ref{sec:const-comm-del} and~\ref{sec:varying-comm-del}, respectively.

	\item[Fusion center:] The central station is in charge of fusing all sensor data to compute a state estimate. We assume that $\mathcal{Z}_t(\tau_p) = \{z_{s_i}^{(i)}(\tau_{p,i}),s_i\in [t_0,t-\tau_{p,i}-\tau_{c,i}]\}_{i\in\mathcal{N}} $ is the dataset available at time $t$ (starting from an initial time $t_0$). Fusion adds further latency, namely the \textit{fusion delay} $\tau_\textit{f,tot}$, which is the sum of 
	the delays $\tau_{f,i}$ required to process the data stream from each sensor $i$.
	In particular, as above,
	we assume that either $\tau_{f,i}$ is constant, or it decreases with the preprocessing delay $\tau_{p,i}$ 
	(intuitively, the more preprocessing is done at the sensor, the less effort is needed for fusion).
	Fig.~\ref{fig:blocks-delays} gives an insight on the processing network with the different latency contributions - by sensor preprocessing, communication, and central station fusion.		
	
\end{description}

\subsection{Optimal Estimation in Processing Networks}

While the sensor data might be received and fused with some (computation and communication) delay, 
we are interested in obtaining an accurate state estimate at the current time $t$; this entails fusing 
sensor information $\mathcal{Z}_t(\tau_p)$ (partially outdated, due to the computation and communication delays) with the open-loop system prediction in~\eqref{eq:processModel}.
This raises a nontrivial \tradeoff:
is it best to transmit raw sensor data and incur larger communication and fusion delays,
or to perform more preprocessing at the edge and transmit more refined (less noisy and more compressed) estimates?
For instance, consider again Fig.~\ref{fig:vehicle-tracking} where robots compute local estimations from images. Each extracted feature both enhances sensor-side accuracy and possibly reduces communication and fusion delays. However, feature extraction comes with preprocessing (edge computation) delay.
A trade-off emerges: on one hand, many features cause a long prediction; on the other hand, few provide a poor estimation.
An \emph{optimal estimation} policy would decide the preprocessing at each sensor in a way to maximize the 
final estimation accuracy.

\textbf{Problem formulation:}
In general, one may wish to optimize the estimation performance at all times, \ie as for Mean Squared Error (MSE) estimators, find $\argmin_{\tau_p \in\Realp{N}} \;\mbox{var}\left(x_{t}-\hat{x}_{t}(\tau_p)\right) $, where $\hat{x}_{t}(\tau_p) \doteq g(\mathcal{Z}_t(\tau_p))$ is a state estimator.
However, such problem comes with the nuisance of time variance. Instead, we resort to its time-invariant steady-state counterpart by exploiting communication reliability (Assumption~\ref{ass:model}).

\begin{prob}\label{prob:time-invariant-P-opt}
	Given the system~\eqref{eq:processModel} with sensor set $\mathcal{N}$ and measurement model~\eqref{eq:measurementModel}, find 
	the optimal preprocessing delays $\tau_p=\{\tau_{p,i}\}_{i\in\mathcal{N}}$ that minimize the steady-state estimation error variance:
	\begin{equation}
	\argmin_{{\small \begin{array}{c}
			\tau_{p,i} \in\mathbb{R}_+, i \in \mathcal{N}
			\end{array}}}  \;\; \psteady(\tau_p)
	\label{prob-2} 
	\end{equation}
	where the total delay is
	\begin{equation}
		\tau_\textit{tot}\doteq \underbrace{\min_{i\in\mathcal{N}}(\tau_{p,i}+\tau_{c,i})}_{\doteq\,\tau_{s}}+\underbrace{\sum_{i\in\mathcal{N}}\tau_{f,i}}_{\doteq\tau_\textit{f,tot}}
		\label{total-delay}
	\end{equation}
	and $\psteady(\tau_p) \doteq \displaystyle\lim_{t\rightarrow +\infty}\!\!\mbox{var}\left(x_{t}-\hat{x}_{t}(\tau_p)\right)$
	is the steady-state estimation error variance. 
	$\tau_\textit{tot}$ accounts for the fact that, due to delays, the steady-state estimate relies on partially outdated measurements: 
	$ \tau_{s} $ is the time it takes to collect data from all sensors (including the freshest available), while $\tau_\textit{f,tot}$ is the time it takes to 
	fuse them.
\end{prob}	

We start by analyzing the single-sensor case to gain some intuition on Problem~\ref{prob:time-invariant-P-opt}.
In the following, to keep notation more readable, we drop the subscripts from the preprocessing delay $\tau_{p,i}$ and refer to it as $\tau$.

\section{Single sensor: preprocessing delay {\large $\tau$}} \label{sec:preprocessing-delay}
	
The goal of this section is twofold: (i) 
to provide a closed-form expression for $\psteady(\tau)$ for the case in which 
we have a single sensor ($N=1$) and we neglect communication and fusion delays ($\tau_\textit{tot} = \tau$), and 
(ii) to compute the optimal preprocessing delay that solves Problem~\ref{prob:time-invariant-P-opt}.
Towards this goal, we use a Kalman filter for state estimation, which is an MSE estimator for linear Gaussian systems.
Moreover, we assume $ \mbox{var}(v_t) \doteq \sigma^2_v(\tau)$ is inversely proportional to preprocessing delay.
This choice is motivated by the observation that the variance of least squares estimation is inversely proportional to the 
number of independent data point (e.g., features extracted in an image). 
Other models are discussed in Section~\ref{sec:other-functions}.

In the single-sensor setup, we can compute the steady-state error variance $\psteadytau(\tau)$ in closed form and 
derive an analytical solution for Problem~\ref{prob:time-invariant-P-opt}.
\begin{thm}[Optimal preprocessing, single sensor]\label{thm-scalar}\hfill\\
	Consider the LTI stochastic scalar system
	\begin{equation}
	\begin{cases}
	dx_t = ax_tdt+dw_t\\
	z_t(\tau) = x_t + v_t(\tau)
	\end{cases} \label{scalar-system}
	\end{equation}
	with state matrix $a \in \mathbb{R}$, process noise $w_t\sim(0,\sigma^2_w)$ with $\sigma^2_w>0$, measurement noise $v_t(\tau)\sim(0,\sigma^2_v(\tau))$ with
	\begin{equation}
	\sigma^2_v(\tau) = \frac{b}{\tau} \qquad b > 0
	\label{sigmav-of-tau}
	\end{equation}
	and initial condition $x_{t_0} \sim (\mu_0,p_{0})$, $ p_0\ge0 $. 
	Assume $\hat{x}_{t}(\tau)$ is the Kalman filter 
	estimate at time $t$ given measurements collected till time $t-\tau$. 
	Then, the steady-state error variance $\psteadytau(\tau)$ has the following expression:
	\begin{equation*}
	\psteadytau(\tau) = \underbrace{\frac{b\mbox{e}^{2a\tau}}{\tau}\left(a+\sqrt{a^2+\sigma^2_w\frac{1}{b}\tau}\right)}_{f(\tau)}+\underbrace{\frac{\sigma^2_w}{2a}\left(\mbox{e}^{2a\tau}-1\right)}_{q(\tau)}
	\end{equation*}
	with limits
	\begin{align}
	\!\begin{split}
	\!\lim_{\tau \rightarrow 0^+} \psteadytau(\tau) \!=\! \lim_{\tau \rightarrow +\infty} \psteadytau(\tau) \!=\! \begin{cases}
	+\infty, &\! a\ge 0\\
	\dfrac{\sigma^2_w}{2|a|}, &\! a<0 \end{cases}
	\end{split}
	\label{error-variance-limits-scalar}
	\end{align} 
	Moreover, $\psteadytau(\tau)$ has a unique global minimum \\ $\tau_\textit{opt}>0$ that satisfies:
	\begin{equation}
	\frac{\sigma^2_w}{b}\tau_\textit{opt}^3 = -a^2\tau_\textit{opt}^2+\frac{1}{4} \label{3rd-degree-eq-scalar-general}
	\end{equation}
\end{thm}
\begin{pf}
	See Appendix \ref{app:proof-thm-scalar}.
\end{pf}

Fig.~\ref{fig:poftautwoparts} illustrates the cost function of Theorem~\ref{thm-scalar}, together with the contributions due to projecting in open-loop the measurement-based estimation and the process noise ($ f(\tau) $ and $ q(\tau) $, respectively).

\begin{figure}
	\centering
	\includegraphics[width=0.7\linewidth]{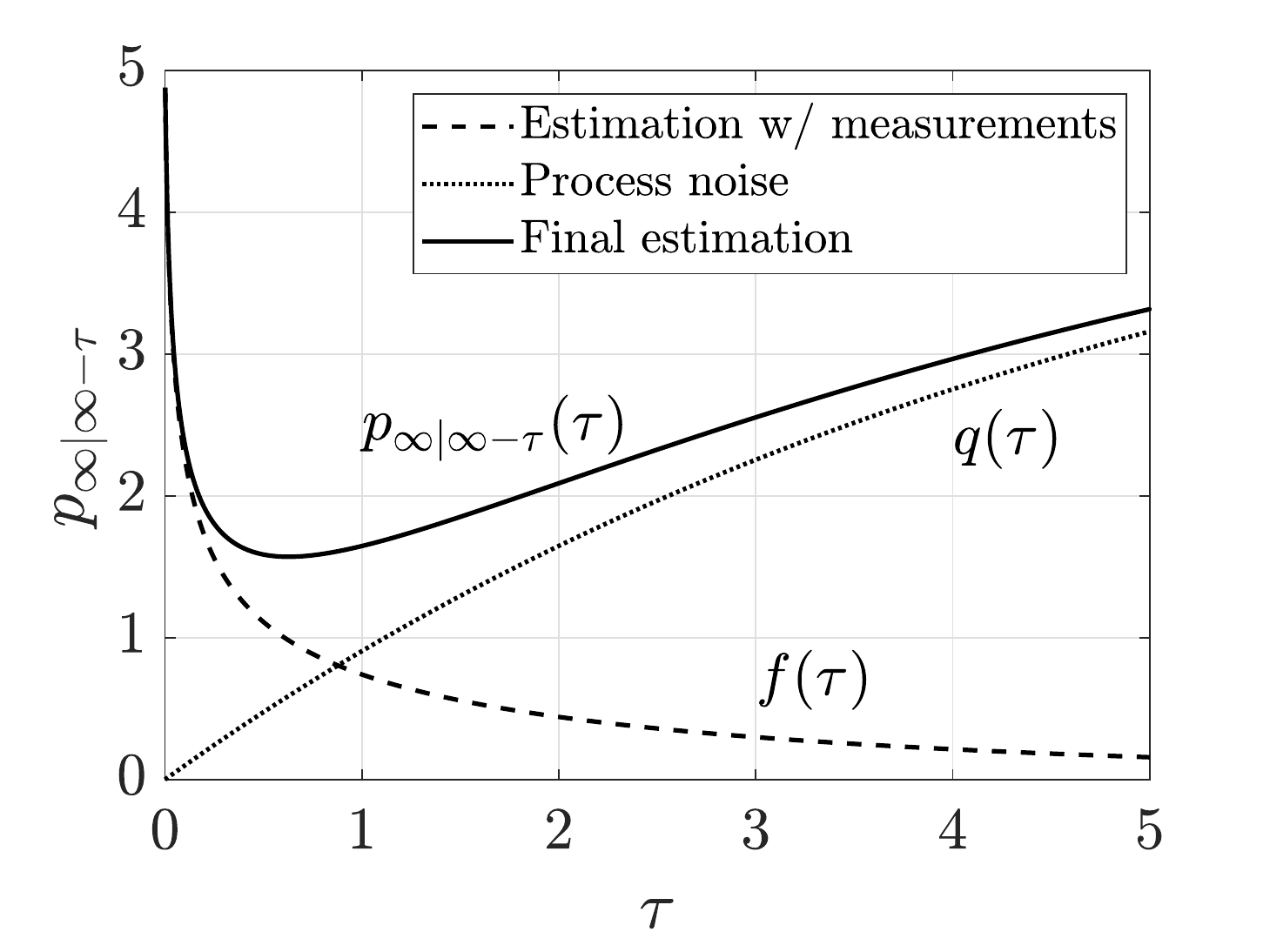}
	\vspace{-3mm}
	\caption{Visual representation of $\psteadytau(\tau)$, with contributions due to estimation $f(\tau)$ and process noise $q(\tau)$.}
	\label{fig:poftautwoparts}
\end{figure}

\subsection{Parameter dependence of optimal delay} \label{sec:param-dependence}

Eq.~\eqref{3rd-degree-eq-scalar-general} provides a characterization for the optimal preprocessing delay $\tau_\textit{opt}$.
Here we discuss how $\tau_\textit{opt}$ behaves as a function of each system's parameter. Notice that $b$ and $ \sigma^2_w $ do not affect $\tau_\textit{opt}$ independently, as they appear in the same coefficient: therefore, we can focus on their ratio $s:=\nicefrac{\sigma^2_w}{b}$. Also, this suggests that what really matters 
is the relative intensity between the process noise and the uncertainty reduction due to preprocessing.

\begin{prop}\label{prop-params}
	$\tau_\textit{opt}$ is strictly decreasing with $s$ and $a^2$. 
\end{prop}

\begin{pf}
	See Appendix~\ref{app:proof-prop-params}.
\end{pf}	
On one hand, Proposition~\ref{prop-params} states that it is more convenient to choose small preprocessing delays for 
``unpredictable systems'', characterized by fast dynamics or large process noise.
On the other hand, if the sensor noise is large, it is convenient to perform further preprocessing,
which explains why $\tau_\textit{opt}$ grows with $b$.

The proof of Proposition~\ref{prop-params} also yields the following upper bound, which may turn useful with uncertain models.
\begin{cor}(Upper bound for $\tau_\textit{opt}$)
	\begin{gather}
	\tau_\textit{opt} \le \tau_{u}\left(|a|,\dfrac{\sigma^2_w}{b}\right) := \begin{cases}
	\frac{1}{2|a|}, & |a| > \sqrt[3]{\dfrac{\sigma^2_w}{2b}}\\
	\sqrt[3]{\dfrac{b}{4\sigma^2_w}}, & \mbox{otherwise}
	\end{cases}\label{tau-upper-bound-scalar-system}
	\end{gather}
\end{cor}

\begin{figure}
	\centering
	\begin{minipage}{.5\textwidth}
		\centering
		\includegraphics[width=.7\textwidth]{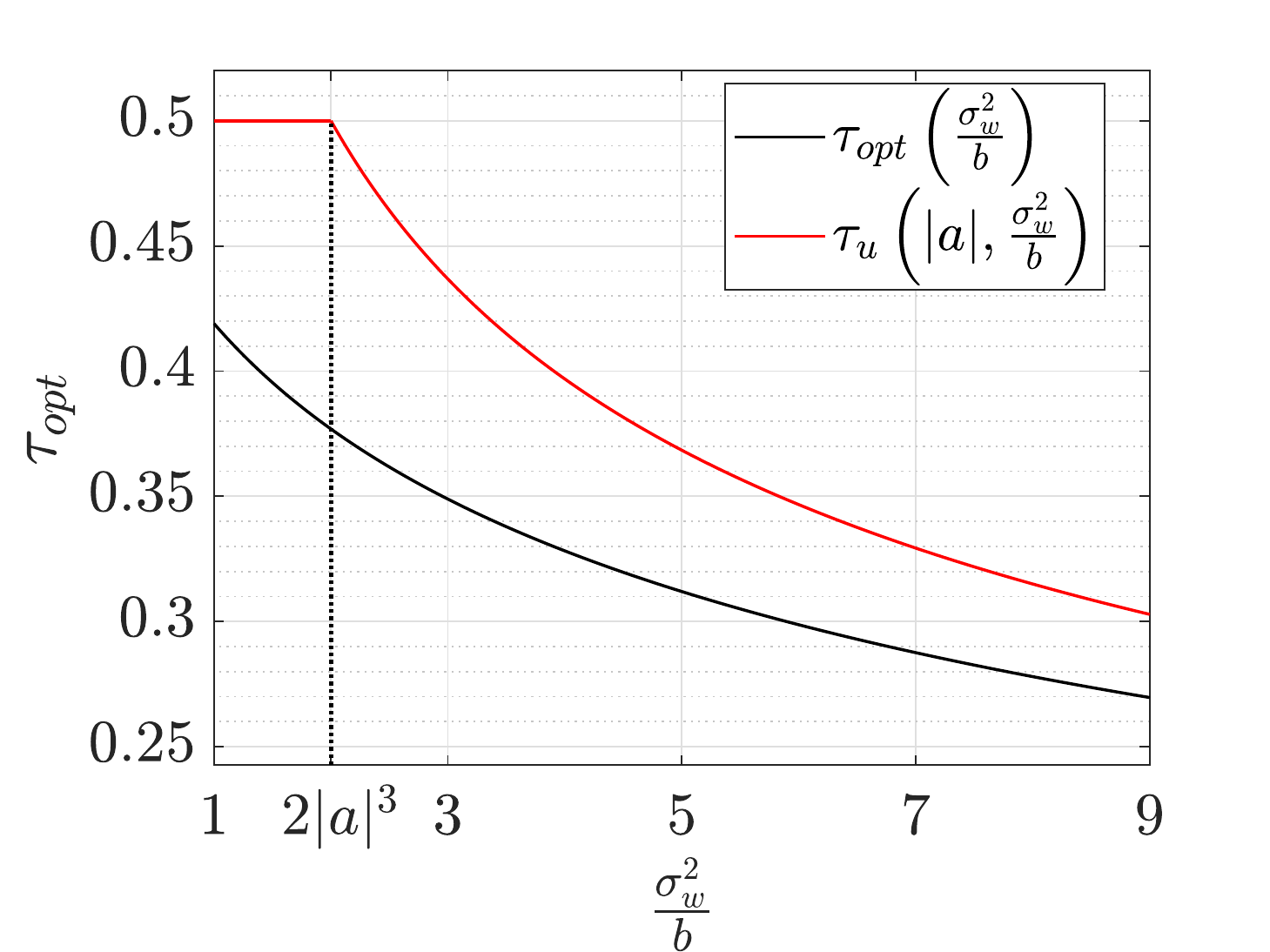}\vspace{-3mm}
		\caption{Optimal delay $\tau_\textit{opt}$ as a function of $ s=\nicefrac{\sigma^2_w}{b} $ \hfill\newline with $a=1$, and upper bound $\tau_{u}$ as per~\eqref{tau-upper-bound-scalar-system}.}
		\label{fig:tau_of_s}
	\end{minipage}
\end{figure}

	Fig.~\ref{fig:tau_of_s} shows how $\tau_\textit{opt}$ varies with $s$ 
	(c.f., eq.~\eqref{tau-of-s-derivative} in Appendix~\ref{app:proof-prop-params}), together with the upper bound in~\eqref{tau-upper-bound-scalar-system}.
	The dependence on $a^2$ is qualitatively similar and omitted for space reasons.

\begin{figure*}
	\centering
	\begin{minipage}{.49\textwidth}
		\centering
		\textbf{Unstable systems}
		\includegraphics[width=.9\linewidth]{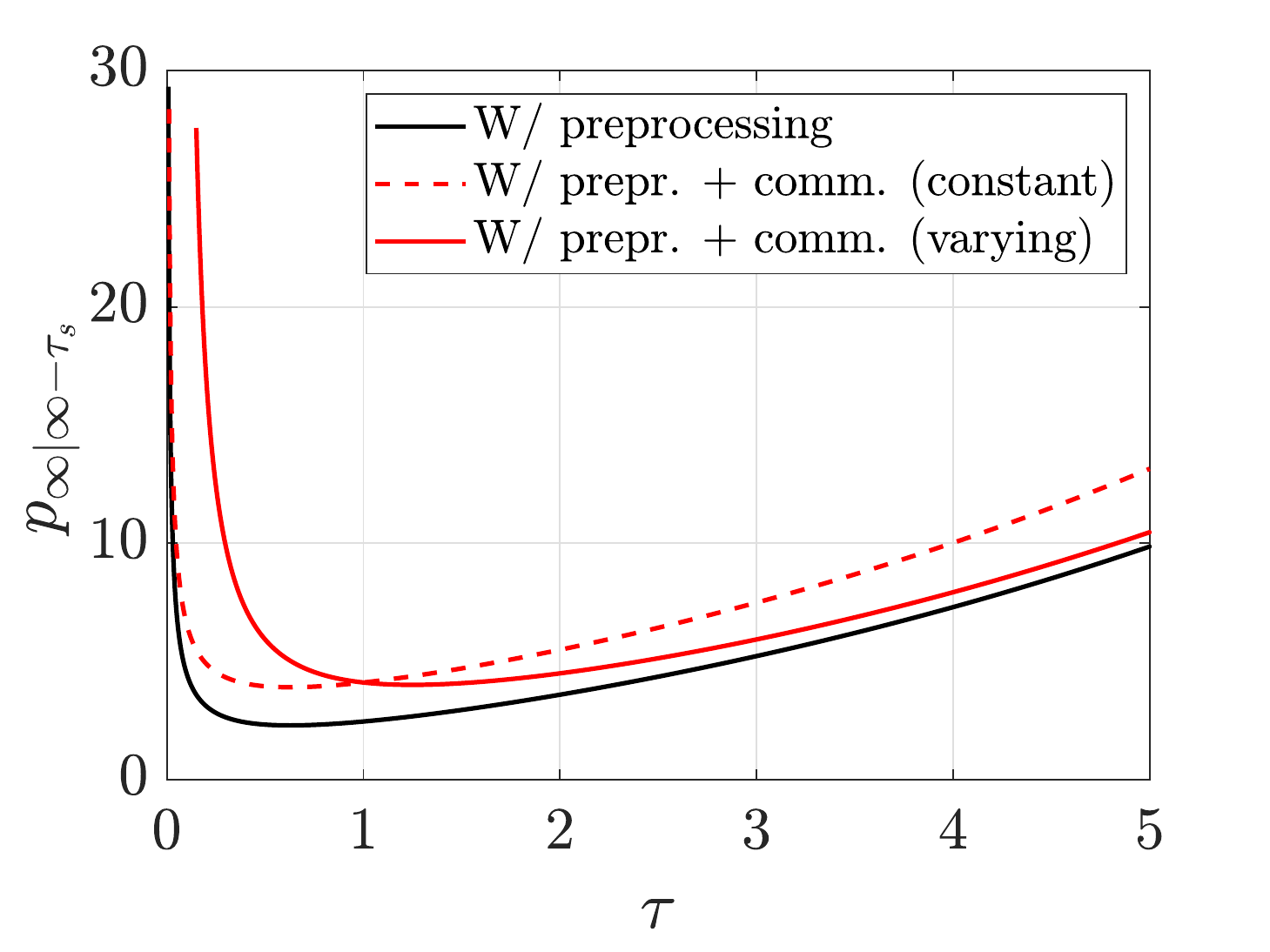}
	\end{minipage}
	\begin{minipage}{.49\textwidth}
		\centering
		\textbf{Stable systems}
		\includegraphics[width=.9\linewidth]{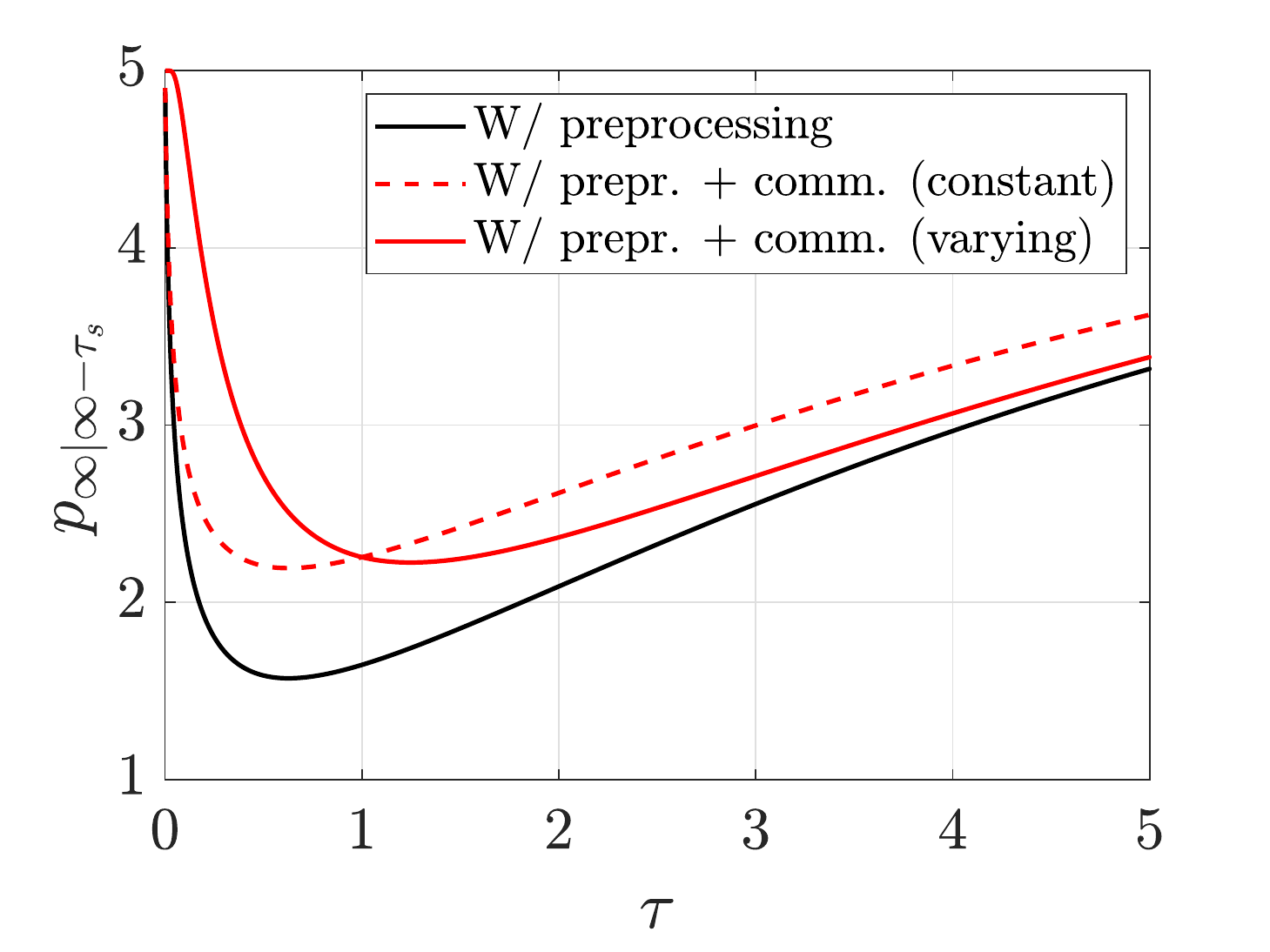}
	\end{minipage}
	\caption{Steady-state error variance $\psteadytautauc(\tau)$ with $\sigma^2_w=b=1$, $a=0.1$ (left) and $a=-0.1$ (right). 
	Black line: no communication delay ($\tau_s=\tau$).
	Dashed red line: constant communication delay ($\tau_c=1$). 
	Solid red line: $\tau$-varying communication delay ($c=1$).}
	\label{fig:p-of-tau-scalar-tot}
\end{figure*}

\subsection{Alternative preprocessing models}\label{sec:other-functions}

Here we consider two different models for the relation between the measurement variance and the preprocessing delay, 
which can be used in place of~\eqref{sigmav-of-tau}.
The models involve a coefficient $\gamma>0$ that can be understood as a convergence rate of an anytime algorithm.
The following case generalizes model~\eqref{sigmav-of-tau} accounting for non-ideality of preprocessing algorithms (as dependent samples).
\begin{cor}[Non-ideal preprocessing]
	Given system (\ref{scalar-system}) and hypotheses as per Theorem~\ref{thm-scalar} with 
	\begin{equation*}
		\sigma^2_v(\tau)=\frac{b}{\tau^{\gamma}} \qquad \gamma> 0
	\end{equation*}
	the steady-state error variance $\psteadytau(\tau)$ has a unique global minimum $\tau_\textit{opt}>0$.
\end{cor}
\begin{pf}
	It can be seen that limits (\ref{error-variance-limits-scalar}) hold and $ \psteadytau(\tau) $ is strictly quasi-convex on $\Realp{}$ (\eg via graphical analysis).
\end{pf}
The second model comes into play with anytime algorithms with exponential convergence, as in~\cite{10.1007/978-3-642-37213-1_16}. 
\begin{cor}[Anytime algorithms]
	Given system (\ref{scalar-system}) and hypotheses as per Theorem~\ref{thm-scalar} with 
	\begin{equation*}
		\sigma^2_v(\tau)=b\mbox{e}^{-\gamma\tau}\qquad \gamma> 0
	\end{equation*}
	the steady-state error variance $\psteadytau(\tau)$ has a unique global minimum $\tau_\textit{opt}>0$ if and only if 
	\begin{equation}
		\gamma > 2\sqrt{\frac{\sigma^2_w}{b}+a^2}  \label{eq:condition-gamma}
	\end{equation}
\end{cor}
\begin{pf}
	In this case, $\tau_\textit{opt}$ can be computed analytically.
	Condition \eqref{eq:condition-gamma} is required to make $\tau_\textit{opt}$ positive.
\end{pf}
\begin{rem}[Phase transition]
	The algorithms whose convergence rate is too slow with respect to the system dynamics are 
	discarded by condition~\eqref{eq:condition-gamma}: if the latter does not hold, $\tau_\textit{opt}=0$, \ie transmitting raw measurements is the optimal choice at sensor side.
\end{rem}

\section{Single sensor: preprocessing and communication delays {\large $\tau,  \tau$}$_{{\tiny c}}$} 
\label{sec:preprocessing-comm-delay}

In this section we add the communication delay, according to the two models mentioned in Section~\ref{sec:set-up}. 
The prediction step therefore stretches to the sensor delay $ \tau_s$ (cf.~\eqref{total-delay}).

\subsection{Constant communication delay}\label{sec:const-comm-del}

In this case, the communication delay $\tau_{c}$ is constant (\ie the preprocessing does not imply data compression): in particular, the transmitted packet number/length is independent on the time spent for preprocessing, which only affects the measurement noise variance. This situation may occur whenever the sensors send quantities whose dimension only depends on the system/algorithms, such as local state estimates.
Being the communication delay constant, it does not impact the optimization with respect to the preprocessing: the steady-state variance $ \psteadytautauc(\tau) $ is simply multiplied by the coefficient $ \mbox{e}^{2a\tau_{c}} $ due to the longer open-loop prediction induced by $\tau_c$. Therefore, the optimal delay is again $\tau_{opt}$ as per Theorem~\ref{thm-scalar}. The dependencies studied in Sec.~\ref{sec:param-dependence} still hold.

\subsection{Computation-dependent communication delay}\label{sec:varying-comm-del}
	
We now turn to the case where the preprocessing also performs data compression,
leading to a $\tau$-varying communication delay which is modelled as:
\begin{equation}\label{comm-delay-comp-dependent}
	\tau_c(\tau) = \dfrac{c}{\tau} \quad c > 0
\end{equation}
with known $ c $.
We have the following result.
\begin{thm}[Optimal preprocessing and communication]\label{thm-scalar-total}
	Given system (\ref{scalar-system}) with measurement noise variance $\sigma^2_v(\tau)$  as per \eqref{sigmav-of-tau},
	and communication delay $\tau_c(\tau)$ as per \eqref{comm-delay-comp-dependent}, 
	the steady-state error variance has expression
	\begin{equation*}
	\psteadytautauc(\tau) = \frac{{b}\mbox{e}^{2a\tau_s }}{\tau}\left(a+\sqrt{a^2+\frac{\sigma^2_w}{{b}}\tau}\right)+\frac{\sigma^2_w}{2a}\left(\mbox{e}^{2a\tau_s }-1\right)\label{p-prep-comm-del}
	\end{equation*}
	 with $ \tau_s=\tau+\nicefrac{c}{\tau} $.
	 Moreover, $\psteadytautauc(\tau)$ admits limits as per~\eqref{error-variance-limits-scalar}, and has a unique global minimum $\tau_{\textit{opt}}>0$.
\end{thm}
\begin{pf}
	See Appendix \ref{app:proof-thm-scalar-total}.
\end{pf}

Fig.~\ref{fig:p-of-tau-scalar-tot} compares the steady-state error variance 
with no communication delay (but with preprocessing delay, in black),
and with 
communication delay (in red, dashed for constant and solid for $\tau$-varying delays) for an unstable and an asymptotically stable systems. Notice that the steepness of the black curve decreasing portion suggests that it is preferable to round $\tau_\textit{opt}$ in excess, if needed, as a lower approximation likely worsens performance. The first communication-delay model (constant $\tau_c$) shifts upward and slightly sharpens the curve, while the second smooths it. In this case, monotonicity of $\tau_\textit{opt}$ as in Section~\ref{sec:param-dependence} cannot be guaranteed. Also, notice that the red curves cross: the model with constant/no compression is outperformed by the $\tau$-varying one if the preprocessing is longer than a minimum threshold.

\begin{rem}
While model~\eqref{comm-delay-comp-dependent} is mainly used for mathematical convenience, in a real setup the 
compression function should be learned or estimated from data.
\end{rem}

\section{Multiple sensors: preprocessing, communication, and fusion delays} \label{sec:multiple-sensors}

We now consider the multi-sensor case and complete our framework by adding the fusion delay.
We model the latter similar to the communication delay: 
$\tau_{f,i}$ is either assumed to be constant, or we assume $\tau_{f,i}(\tau_i)=\nicefrac{f_i}{\tau_i}$, where $f_i>0$ is a known constant.
We consider for simplicity $N$ identical independent sensors, each with the same delays $\tau$ (preprocessing), $ \tau_c(\tau) $ (communication) and $\tau_{f}(\tau)$ (fusion), the latter two being constant or varying.
The overall delay~\eqref{total-delay} becomes 
\begin{equation}\label{eq:totDelayHom}
	\tau_\textit{tot} = \underbrace{\tau+\tau_{c}(\tau)}_{\tau_{s}}+\underbrace{\tau_{f}(\tau)N}_{\tau_\textit{f,tot}}
\end{equation}
and the network measurements in~\eqref{eq:measurementModel} are then modeled as
\begin{equation}
z_{t}(\tau) = \begin{bmatrix}
1 &
\dots & 
1
\end{bmatrix}^Tx_{t} + v_{t}(\tau) \qquad R(\tau) = \frac{b}{\tau}I_N
\label{n-sensors}
\end{equation}

From the least squares framework, it is well known that such system yields an overall variance reduction for $ z_t(\tau) $ which is linear with the number of samples.
Alternatively, the homogeneous $ N $-sensor network can be seen (from the standpoint of the estimation performance) as a single sensor with processing noise variance $\sigma^2_v(\tau)$ reduced by a factor $ N $ with respect to each sensor in~\eqref{n-sensors}, and total delay~\eqref{eq:totDelayHom}. Then, the optimal computational delay for such virtual single sensor also maximizes the performance for the network~\eqref{n-sensors}.

\begin{rem}
	The advantage of multi-sensor networks is reducing the measurement-noise variance for each sensor, yielding more accurate state estimation.
\end{rem}

\begin{rem}
A common wisdom in estimation theory is that adding more sensors always yields better estimates (possibly with performance saturation). We show that, if fusion time has to be considered, the optimal solution is adding sensors up to a certain amount, since the cost of processing more overtakes the sample variance reduction.
\end{rem}

\begin{figure}
	\centering
	\begin{tikzpicture}[scale=.9]
	\begin{axis}[grid=major,domain=1:10,ymax=3.4,ylabel=$\psteady$,xlabel=Number of sensors $ S $,xtick={1,2,3,4,5,6,7,8,9,10},legend style={at={(.2,.9)},anchor=west}]
	\addplot+[mark=square*,mark size=3,mark options={fill=black},samples=10,color=black] {1/10*exp(-2*.1-2*.1-2*0.02*x)*(1/.1/x)*(-1+sqrt(1+x*10*.1*10))-10/2*(exp(-2*.1-2*.1-2*0.02*x)-1)};
	\addlegendentry{W/ fusion delay}
	\addplot+[mark=square*,mark size=3,mark options={fill=red},samples=10,color=red] {1/10*exp(-2*.1-2*.1)*(1/.1/x)*(-1+sqrt(1+x*10*.1*10))-10/2*(exp(-2*.1-2*.1)-1)};
	\addlegendentry{W/o fusion delau}
	\end{axis}
	\end{tikzpicture}
	\caption{Variance $\psteady(S)$ with $a=-1$, $\sigma^2_w=10$, $b=\tau=0.1$, $\tau_c=0.1$, $\tau_f=0.02$ (black) and  $\tau_f=0$ (red, no fusion delay).}
	\label{fig:p-of-n-scalar}
\end{figure}
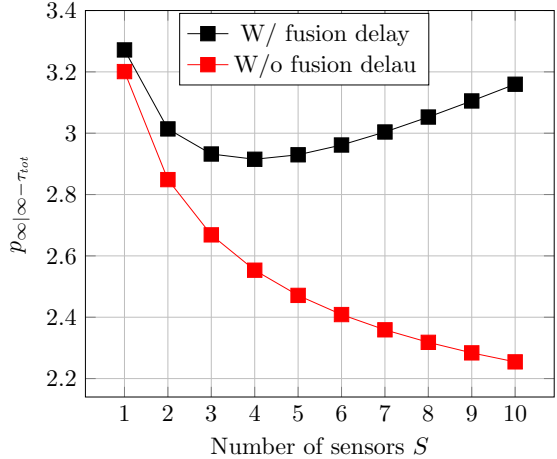

If $ \psteady $ is seen as a discrete function of the number of sensors $ S $ ($S\le N$) with fixed $\tau$, 
one can also ask if there is an \emph{optimal sensor quantity} $ S_{opt} $, corresponding to $\argmin \psteady(S) $\footnote{Due to its discrete domain, $ \psteady(S) $ may have two points of global minimum with $\psteady(S^*)=\psteady(S^*+1)$.}.
Fig.~\ref{fig:p-of-n-scalar} shows the performance behaviour according to the number of sensors, with fixed $\tau$. Notice that neglecting the fusion delay may yield important performance drops (about 12\% with $ S_{opt} $ sensors and 32\% with all $ N $). 
Then, Problem~\ref{prob:time-invariant-P-opt} can be extended
to decide on the optimal number of sensors to be used for estimation.

\begin{prob}[Homogeneous network]\label{prob:hom-network}
	Given system (\ref{scalar-system}) with $N$ identical sensors and measurement model~\eqref{n-sensors}, find the optimal sensor amount S and preprocessing delay $\tau$ that minimize the steady-state error variance:
	\begin{equation*}
	\argmin_{{\small \begin{array}{c}
		S\in\{1,...,N\}\\
		\tau\in\mathbb{R}_+
		\end{array}} }\psteady(S,\tau)
	\end{equation*}
\end{prob}

\begin{figure}
	\centering
	\begin{tikzpicture}[scale=.9]
	\begin{axis}[grid=major,domain=1:10,ymax=3.6,ylabel=$\psteady$,xlabel=Number of sensors $ S $,xtick={1,2,3,4,5,6,7,8,9,10},legend style={at={(.3,.8)},anchor=west}]
	\addplot+[mark=square*,mark size=3,mark options={fill=red},samples=10,color=red] {1/10*exp(-2*(.05+.1+0.02*x))*(1/.05/x)*(-1+sqrt(1+x*10*.05*10))-10/2*(exp(-2*(.05+.1+0.02*x))-1)};
	\addlegendentry{$\tau=0.05$}
	\addplot+[mark=square*,mark size=3,mark options={fill=black},samples=10,color=black] {1/10*exp(-2*(.1+.1+0.02*x))*(1/.1/x)*(-1+sqrt(1+x*10*.1*10))-10/2*(exp(-2*(.1+.1+0.02*x))-1)};
	\addlegendentry{$\tau=0.1$}
	\addplot+[mark=square*,mark size=3,mark options={fill=green},samples=10,color=green] {1/10*exp(-2*(.15+.1+0.02*x))*(1/.15/x)*(-1+sqrt(1+x*10*.15*10))-10/2*(exp(-2*(.15+.1+0.02*x))-1)};
	\addlegendentry{$\tau=0.15$}
	\addplot+[mark=square*,mark size=3,mark options={fill=blue},samples=10,color=blue]
	{1/10*exp(-2*(.2+.1+0.02*x))*(1/.2/x)*(-1+sqrt(1+x*10*.2*10))-10/2*(exp(-2*(.2+.1+0.02*x))-1)};
	\addlegendentry{$\tau=0.2$}
	\end{axis}
	\end{tikzpicture}
	\caption{Variance $\psteady(S)$ with $a=-1$, $\sigma^2_w=10$, $b=\tau_c=0.1$, $\tau_f=0.02$, $ \tau\in\{0.05,0.1,0.15,0.2\} $.}
	\label{fig:p-of-n-tau-scalar-const}
\end{figure}
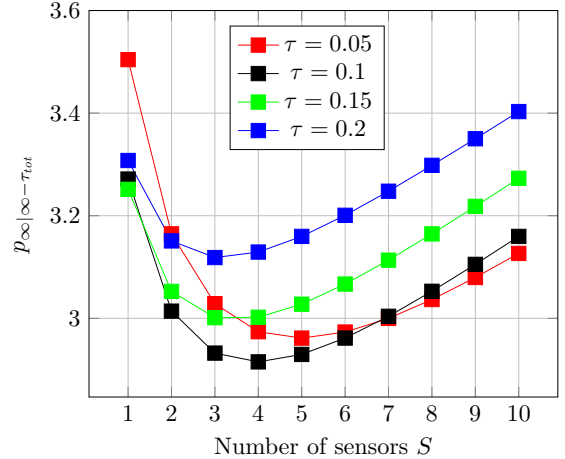

Proving uniqueness of the solution in this case in nontrivial, due to both the discrete domain of the cost function ($ S $ must be natural) and the difficulty of proving quasi-convexity (or a suitably equivalent characterization).
However, simulations results suggest that Problem~\ref{prob:hom-network} admits a unique solution:
Fig.~\ref{fig:p-of-n-tau-scalar-const} shows $ \psteady(S) $ corresponding to different values of $\tau$ with constant delay. The $\tau$-varying model is similar and omitted for space reasons.

\section{Extensions to heterogeneous, multi-variate, and discrete-time systems}\label{sec:discussion}

This section presents extensions to heterogeneous networks and discrete-time systems, as well as future work directions. The interested reader can find more details in~\cite{2019arXiv191105859B}, where some realistic scenarios are analyzed.

\subsection{Sensor selection in heterogeneous networks}\label{sec:het-nets}

In general, sensors in the processing network might have different resources, resulting in different coefficients $ b $ in~\eqref{sigmav-of-tau}. 
If this is the case, computing the optimal preprocessing becomes a multivariate problem.
Moreover, if sensors are heterogeneous, one also faces the choice of whether to use the data from all sensors or
disregard data from some of them. Therefore, a potential generalization of Problems~\ref{prob:time-invariant-P-opt} and 
Problem~\ref{prob:hom-network} is as follows.

\begin{prob}[Heterogeneous network]\label{prob-multisensor}
	Given system (\ref{scalar-system}) with sensor set $ \mathcal{N} $ and measurement model~\eqref{eq:measurementModel}, find the optimal sensor subset $\mathcal{S}$ and preprocessing delays $ \tau=\{\tau_i\}_{i\in\mathcal{S}} $ that minimize the steady-state error variance:
	\begin{equation*}
	\argmin_{{\small \begin{array}{c}
		\mathcal{S}\subseteq\mathcal{N}\\
		\tau_i\in\mathbb{R}_+,i\in\mathcal{S}
		\end{array}}}\psteady(\mathcal{S},\tau)
	\end{equation*}
\end{prob}
Note that the combinatorial nature of the problem makes it difficult to compute exact solutions. In the extended paper~\cite{2019arXiv191105859B} we investigate this formulation with simulations, and propose approximate algorithms.

\subsection{Discrete-time and multi-dimensional systems}\label{sec:multidim}

While the continuous-time framework was instrumental to obtain insights and analytical solutions, in practical problems it is interesting to consider a discrete-time formulation due to the digital nature of involved systems and algorithms. In~\cite{2019arXiv191105859B}, we extend the setup considered in this paper to more general scenarios, accounting for discrete-time and multi-dimensional states.
Numerical simulations confirm that the trends observed in the scalar case also arise in such more general case: 
the \tradeoff can be optimized by suitably selecting sensors and preprocessing delays.

\subsection{Dealing with channel constraints}\label{sec:comm-constraints}

An interesting avenue for future research is to consider more realistic communication model, 
including 
finite bandwidth, channel capacity, unreliability, or packet loss. For instance, to model limited channel capacity, an upper bound may be imposed on the total communication delay:
\begin{equation*}
\sum_{i\in\mathcal{S}} \tau_{c,i}(\tau_i) \le {\tau}_{u,c}
\end{equation*}
In this way, each sensor is forced not to keep the channel busy for too long, letting all sensors transmit their data.

\section{Conclusions}\label{sec:conclusions}

In this paper, we investigate optimal estimation in a processing network in the 
presence of communication and computational delays.
We model sensor-side preprocessing as a stochastic measurement model, 
whose noise intensity decreases with the computational delay.
Similarly, communication and fusion delays are modeled as a constant or decreasing function of computation delay, simulating data compression.
For the continuous-time, scalar, single-sensor scenario, we prove that the resulting trade-off between
preprocessing and computation can be optimized analytically.
We further extend these results to the case of a network of homogeneous sensors, where one has 
also to account for the fusion delay incurred at the central station which is in charge of fusing all the sensor measurements.
We conclude the paper by discussing several ongoing efforts to extend this work to the case of a 
multi-variate, heterogeneous processing networks, monitoring a discrete-time system.


	\appendix

	\section{Proof of Theorem \ref{thm-scalar}}\label{app:proof-thm-scalar}
	
	By considering model~\eqref{sigmav-of-tau}, the steady-state Kalman error variance associated with $ \hat{x}_{t-\tau}(\tau) $ (outdated estimate) is
	\begin{equation}
	p_\infty(\tau) = \frac{{b}}{\tau}\left(a+\sqrt{a^2+\frac{\sigma^2_w}{{b}}\tau}\right) \label{p-infty}
	\end{equation}	
	The model-based open-loop predictor error has dynamics
	\begin{equation}
	d\Tilde{x}_{s}(\tau) = a\Tilde{x}_{s}(\tau)ds + dw_s, \quad t-\tau\le s\le t
	\label{open-loop-error-dynamics-scalar}
	\end{equation}
	Then, the error at time $t$ is given by solving (\ref{open-loop-error-dynamics-scalar}) as a Cauchy problem with initial condition $\Tilde{x}_{t-\tau}(\tau)$:
	\begin{equation*}
	\Tilde{x}_{t}(\tau) = \mbox{e}^{a\tau}\Tilde{x}_{t-\tau}(\tau)+\Bar{w}(\tau)
	\end{equation*}
	where $\Bar{w}(\tau)$ is the stochastic integral of $ w_s $ in $[t-\tau,\: t]$. The steady-state prediction error variance is then 
	\begin{align*}
	\begin{split}
	\psteadytau(\tau) &\overset{(i)}{=} \mbox{var}(\mbox{e}^{a\tau}\Tilde{x}_{t-\tau}) + \mbox{var}(\Bar{w}(\tau)) =\\
	& = \frac{{b}\mbox{e}^{2a\tau}}{\tau}\left(a+\sqrt{a^2+\frac{\sigma^2_w}{{b}}\tau}\right)+\frac{\sigma^2_w}{2a}\left(\mbox{e}^{2a\tau}-1\right)
	\end{split}
	\end{align*}
	where $(i)$ is motivated by uncorrelated terms. Indeed, $\tilde{x}_{t-\tau} \in \mbox{span}\{x_{t_0}, w_s, v_s : t_0 \le s\le t-\tau\}$, while $\Bar{w}(\tau) \in \mbox{span}\{w_s, \: t-\tau \le s \le t\}$, with $w_t$ white noise and $w_t \perp x_{t_0}, v_s \ \forall t\ge t_0,\forall s$ by hypothesis. The only sample providing nonzero correlation is $w_{t-\tau}$, but having zero Lebesgue measure its contribution to $\Bar{w}_t$ is null.\\
	We proceed now in studying critical points of $ \psteadytau(\tau) $, since being limits~\eqref{error-variance-limits-scalar} equal at both domain extrema and being $  \psteadytau(\tau)  \in\mathcal{C}^0(\bar{\mathbb{R}}_+)$ at least one must exist. By setting $ \psteadytau'(\tau)=0 $ and rejecting $\tau=0$, we get
	\begin{equation}
	\frac{\sigma^2_w}{b}\tau^3 + a^2\tau^2-\frac{1}{4} = 0 
	\label{3rd-degree-eq-scalar-general-1}
	\end{equation}
	
	This equation always admits a real positive solution, in virtue of Bolzano's theorem by considering $F(\tau):=\nicefrac{\sigma^2_w}{b}\tau^3 + a^2\tau^2-\nicefrac{1}{4}$ and $I:=\left[0,\nicefrac{1}{2|a|}\right]$. Strict quasi-convexity of $ \psteadytau(\tau) $ on $\Realp{}$ can be checked via convexity of its sublevel sets, \eg with a graphical analysis. Such property guarantees that the critical point is indeed the unique point of global minimum.
		\section{Proof of Proposition \ref{prop-params}} \label{app:proof-prop-params}
	
	For this proof we are going to exploit the implicit function theorem, whose statement is recalled for convenience.
	\begin{thm}[Dini's theorem]
		Let F be a continuously differentiable function on some open $D \subset \mathbb{R}^2$. Assume that there exists a point $(\bar{x},\bar{y}) \in D$ such that:
		\begin{itemize}
			\item $F(\bar{x},\bar{y})=0$;
			\item $\dfrac{\partial F}{\partial y}(\bar{x},\bar{y})\neq 0$.
		\end{itemize}
		Then, there exist two positive constant a, b and a function $f : I_{\bar{x}}:=(\bar{x}-a,\bar{x}+a) \mapsto J_{\bar{y}}:=(\bar{y}-b,\bar{y}+b)$ such that
		\[F(x,y)=0 \iff y = f(x) \quad \forall x \in I_{\bar{x}}, \ \forall y \in J_{\bar{y}}\]
		Moreover, $f \in \mathcal{C}^1(I_{\bar{x}})$ and
		\begin{equation}\label{Dini-thm-der}
			f'(x) = -\frac{F_x(x,f(x))}{F_y(x,f(x))} \quad \forall x \in I_{\bar{x}}
		\end{equation}
		\label{thm-Dini}
	\end{thm}
	We can see the left-hand term in eq. (\ref{3rd-degree-eq-scalar-general}) as a one-parameter function of two positive-valued variables, namely
	\[ F: \mathbb{R}_+ \times \mathbb{R}_+ \rightarrow \mathbb{R}, \ (\pi,\tau) \mapsto F(\pi,\tau) = s\tau^3 + a^2\tau^2-\frac{1}{4}\]
	with $\pi=\{s,a^2\}$. First of all, let us check if Dini's theorem hypotheses are satisfied: given a solution $(\bar{\pi}, \bar{\tau}_{opt})$ of eq. (\ref{3rd-degree-eq-scalar-general}) it holds
	\begin{itemize}
		\item $F(\bar{\pi}, \bar{\tau}_{opt}) = 0$, by construction;
		\item $F_{\tau}(\bar{\pi}, \bar{\tau}_{opt}) = 3s\bar{\tau}_{opt}^2 +2a^2\bar{\tau}_{opt} > 0$, since all variables are positive.
	\end{itemize}
	Then Theorem~\ref{thm-Dini} applies and there exists a function $\tau(\pi)$ such that $F(\pi,\tau_{opt}) = 0 \iff \tau_{opt} = \tau(\pi)$, with $\pi$ in some open neighbourhood of $\bar{\pi}$. In fact, since we did not pose constraints on $(\bar{\pi}, \bar{\tau}_{opt})$, we can state that such a function is defined for all $ \pi\in\mathbb{R}_+ $. The two cases for $\pi$ are studied independently.
	\begin{description}
		\item[$\boldsymbol{\pi = s}$] 
		By~\eqref{Dini-thm-der}, the first derivative of $\tau(\pi) = \tau(s)$ is
		\begin{equation}
		\tau'(s) = -\dfrac{F_{s}(s,\tau(s))}{F_{\tau}(s,\tau(s))} = -\dfrac{{\tau(s)}^2}{3{s}{\tau(s)} +2a^2} \label{tau-of-s-derivative}
		\end{equation}
		We conclude that $\tau'(s) < 0 \ \forall s \in \mathbb{R}_+$, namely, $\tau_{opt}$ is strictly decreasing with $s$.
		
		\item[$\boldsymbol{\pi = a^2}$] The first derivative of $\tau(\pi) = \tau(a^2)$ is
		\begin{align}
		\tau'(a^2) &= -\dfrac{F_{a^2}(a^2,\tau(a^2))}{F_{\tau}(a^2,\tau(a^2))} = -\dfrac{\tau(a^2)}{3s{\tau(a^2)} +2a^2} \label{tau-of-a-derivative}
		\end{align}
		We conclude that $\tau'(a^2) < 0 \ \forall a \in \mathbb{R}$, namely, $\tau_{opt}$ is strictly decreasing with $a^2$ (regardless of sign$(a)$).
	\end{description}
	\section{Proof of Theorem \ref{thm-scalar-total}} \label{app:proof-thm-scalar-total}
	
By applying open-loop prediction to $ p_{\infty}(\tau) $~\eqref{p-infty} so as to cover the delay $ \tau_s=\tau+\tau_c $, the steady-state prediction error variance takes the expression in Theorem~\ref{thm-scalar-total}. Strict quasi-convexity holds also in this case (again, this can be shown with graphical analysis). In virtue of this fact, continuity and limits \eqref{error-variance-limits-scalar}, we conclude that the point of minimum exists unique and is different from zero.


\end{document}